\documentclass[12pt]{amsart}
\usepackage{latexsym,amsfonts,amsthm,amsmath,amscd,amssymb}
\usepackage[dvips]{graphicx}
\usepackage{graphicx}
\usepackage{srcltx}

\reversemarginpar

\newcommand{\matN} {{\mathbb N}}
\newcommand{\matZ} {{\mathbb Z}}
\newcommand{\matR} {{\mathbb R}}

\newcommand{\matP} {{\mathbb P}}

\newcommand{\matH} {{\mathbb H}}
\newcommand{\vol} {\mathrm{vol}}
\newcommand{\finedimo}{{\hfill\hbox{$\square$}\vspace{2pt}}}

\newcommand{\myparagraph}[1]{\medskip\noindent\textbf{#1}\quad}

\newtheorem{lemma}{Lemma}[section]
\newtheorem{thm}[lemma]{Theorem}
\newtheorem{rem}[lemma]{Remark}
\newtheorem{prop}[lemma]{Proposition}

\begin{document}

\author[Petronio]{Carlo Petronio}
\address{\noindent Dipartimento di Matematica Applicata\\ Universit\`a di Pisa\\
Via Fi\-lip\-po Buonarroti 1C\\ 56127, Pisa, Italy}
\email{petronio@dm.unipi.it}

\author[Vesnin]{Andrei Vesnin}
\address{Sobolev Institute of Mathematics\\ Novosibirsk 630090\\
Russia} \email{vesnin@math.nsc.ru}

\address{Department of Mathematics\\ Omsk Technical University\\ pr. Mira 11\\ Omsk 644050\\ Russia} \email{vesnin@omgtu.ru}

\title[Bounds for the complexity of branched coverings]
{Two-sided bounds for the complexity\\ of cyclic branched
coverings\\ of two-bridge links}

\thanks{This work is the result of a collaboration
carried out in the framework of the INTAS project ``CalcoMet-GT''
03-51-3663. The second author was also partially supported by the grant
NSh-8526.2006.1, by a grant of the RFBR, and by a grant of the Siberian Branch
of RAN}


\begin{abstract}
We consider closed orientable 3-dimensional hyperbolic manifolds
which are cyclic branched coverings of the 3-sphere, with branching
set being a two-bridge knot (or link).
We establish two-sided linear bounds depending on the order of the covering
for the Matveev complexity of the covering ma\-ni\-fold. The lower
estimate uses the hyperbolic volume and results of Cao-Meyerhoff and
Gu\'{e}ritaud-Futer (who recently improved previous work of
Lackenby), while the upper estimate is based on an explicit
triangulation, which also allows us to give a bound on the Delzant
T-invariant of the fundamental group of the manifold.
\end{abstract}

\date{\today}

\maketitle

\section{Definitions, motivations and statements}
\myparagraph{Complexity} Using simple spines (a technical notion
from piecewise linear topology that we will not need to recall in
this paper), Matveev~\cite{Acta} introduced a notion of
\emph{complexity} for compact 3-dimen\-sio\-nal manifolds. If $M$ is
such an object, its complexity $c(M)\in\matN$ is a very efficient
measure of ``how complicated'' $M$ is, because:
\begin{itemize}
\item every 3-manifold can be uniquely expressed as a connected
sum of prime ones (this is an old and well-known
fact, see~\cite{hempel});
\item $c$ is additive under connected sum;
\item if $M$ is closed and prime, $c(M)$ is precisely the minimal number of
tetrahedra needed to triangulate $M$.
\end{itemize}
In the last item the notion of triangulation is only meant in a
loose sense, namely just as a gluing of tetrahedra along faces, and
an exception has to be made for the four prime $M$'s for which
$c(M)=0$, that is $S^3$, $\matR\matP^3$, $S^2\times S^1$,
and $L(3,1)$.

Computing exactly the complexity $c(M)$ of any given 3-manifold $M$
is theoretically very difficult, even if quite easy experimentally,
using computers~\cite{DAN05}. In the closed prime case the state of
the art is as follows:
\begin{itemize}
\item A computer-aided tabulation of the closed $M$'s with
$c(M)\leqslant 12$ has been completed in various
steps~\cite{MPexp,DAN05,c12} (see also~\cite{matbook});
\item A general lower bound for $c(M)$ in terms of the homology of
$M$ was established in~\cite{mape};
\item Asymptotic two-sided bounds for the complexity of some
specific infinite series of manifolds were obtained
in~\cite{MaPV,pepe};
\item A conjectural formula for the complexity of any Seifert
fibred space and torus bundle over the circle
was proposed (and proved to be an upper bound)
in~\cite{MPded}.
\end{itemize}

Several other results, including exact computations for infinite
series, have been obtained in the case of manifolds with non-empty
boundary, see~\cite{anis,fmp1,fmp2,fmp3}. Since we will stick in
this paper to the closed case, we do not review them here.

Using the hyperbolic volume and deep results of Lackenby~\cite{L}
improved recently for the case of hyperbolic two-bridge links
in~\cite{GF}, and of Cao-Meyerhoff~\cite{CM}, together with explicit
triangulation methods to be found~\cite{Mi,MuV}, we will analyze in
this paper the complexity of cyclic coverings of the 3-sphere
branched along two-bridge knots and links. More specifically, we
will prove asymptotic
two-sided linear estimates for the complexity in terms of the order
of the covering. Before giving our statements we need to recall some
terminology.

\myparagraph{Two-bridge knots and links} If $p,q$ are coprime
integers with $p\geqslant 2$ we denote by $K(p,q)$ the two-bridge
link in the 3-sphere $S^3$ determined by $p$ and $q$,
see~\cite{BZ,Kawa,MuV}. It is well-known that $K(p,q)$ does not
change if a multiple of $2p$ is added to $q$, so one can assume
that $|q|<p$. In addition $K(p,-q)$ is the mirror image of
$K(p,q)$. Therefore, since we will not care in the sequel about
orientation, we can assume $q>0$. Summing up, from now on our
assumption will always be that the following happens:
\begin{equation}\label{pq:hp}
p, q \in \mathbb Z, \qquad p \geqslant  2, \qquad 0 < q < p, \qquad
(p,q)=1.
\end{equation}

We recall that if $p$ is odd then $K(p,q)$ is a knot, otherwise it
is a 2-component link; moreover, two-bridge knots and links are
alternating~\cite[p.~189]{BZ}. Planar alternating diagrams of
$K(p,q)$ will be shown below. Since we are only interested in the
topology of the branched coverings of $K(p,q)$, we regard it as an
unoriented knot (or link), and we define it to be equivalent to some
other $K(p',q')$ if there is an automorphism of $S^3$, possibly an
orientation-reversing one, that maps $K(p,q)$ to $K(p',q')$. It is
well-known (see~\cite[p.~185]{BZ}) that $K(p',q')$ and $K(p,q)$ are
equivalent if and only if $p'=p$ and $q'\equiv \pm q^{\pm1}$ (mod
$p$).

Under the current assumption~(\ref{pq:hp}), the two-bridge knot (or
link) $K(p,q)$ is a torus knot (or link) precisely when $q$ is $1$
or $p-1$, and it is hyperbolic otherwise. The simplest
non-hyperbolic examples are the Hopf link $K(2,1)$, the left-handed
trefoil knot $K(3,1)$ and its mirror image $K(3,2)$, the
right-handed trefoil (but we are considering a knot to be equivalent
to its mirror image, as just explained). The easiest hyperbolic
$K(p,q)$ is the figure-eight knot $K(5,2)$.

\myparagraph{Branched coverings} If $K(p,q)$ is a knot
(\emph{i.e.}~$p$ is odd) and $n\geqslant 2$ is an integer, the
$n$-fold cyclic covering of $S^3$ branched along $K(p,q)$ is a
well-defined closed orientable 3-manifold that we will denote by
$M_n(p,q)$. One way of defining it is as the metric completion of
the quotient of the universal covering of $S^3\setminus K(p,q)$
under the action of the kernel of the homomorphism
$\pi_1(S^3\setminus K(p,q))\to\matZ/_{n\matZ}$ which factors
through the Abelianization $\pi_1(S^3\setminus K(p,q))\to
H_1(S^3\setminus K(p,q))$ and sends a meridian of $K(p,q)$, which
generates $H_1(S^3\setminus K(p,q))$, to $[1]\in\matZ/_{n\matZ}$.

If $K(p,q)$ is a link and a generator $[m]$ of $\matZ/_{n\matZ}$ is
given, a similar construction defines the \emph{meridian-cyclic}
branched covering $M_{n,m}(p,q)$ of $S^3$ along $K(p,q)$, by
requiring the meridians of the two components of $K(p,q)$ to be sent
to $[1]$ and $[m]\in\matZ/_{n\matZ}$ respectively. Note that
meridian-cyclic coverings are also called \emph{strongly cyclic}
in~\cite{Zim}, and that the two components of $K(p,q)$ can be
switched, therefore we do not need to specify which meridian is
mapped to $[1]$ and which to $[m]$. Since in the sequel we will
prove estimates on the complexity of $M_{n,m}(p,q)$ which depend on
$n$ only and apply to every $M_{n,m}(p,q)$, with a slight abuse we
will simplify the notation and indicate by $M_n(p,q)$ an arbitrary
meridian-cyclic $n$-fold covering of $S^3$ branched along $K(p,q)$.
This will allow us to give a unified statement for knots and links.
We recall that $M_2(p,q)$ is the lens space $L(p,q)$.

\myparagraph{Continued fractions}
In the sequel we will employ continued fractions, that we define as follows:
$$[a_1,a_2,\ldots,a_{k-1},a_k]=a_1 +
\frac{\displaystyle 1}{\displaystyle a_2 + \, \cdots \, +
\frac{\displaystyle 1}{\displaystyle a_{k-1} + \frac{\displaystyle
1}{\displaystyle a_{k}}}}.$$
Given $p,q$ satisfying~(\ref{pq:hp}), we now recall~\cite[p.~25]{Kawa}
that there is a unique
minimized expansion of $p/q$ as a continued fraction with positive entries,
namely an expression as
$p / q = [a_1,\ldots,a_k]$
with $a_1,\ldots,a_{k-1}>0$ and $a_k>1$. (The expansion is called \emph{minimized}
because if $a_k=1$ then $[a_1,\ldots,a_{k-1},1]=[a_1,\ldots,a_{k-1}+1]$,
as one easily sees.)
We then define $\ell(p,q)$ to be $k$ if $a_1>1$ and $k-1$ if $a_1=1$.

This apparently original definition of $\ell(p,q)$ is explained by
the following result established below (see also the proof of
Proposition~\ref{twist:prop}):

\begin{prop}\label{strange:ell:prop}
$\ell(p,q)$ is the minimum of the
lengths of positive continued fraction expansions
of rational numbers $p'/q'$ such that $K(p',q')$ is equivalent to $K(p,q)$.
\end{prop}

\begin{rem}
$\ell(p/q)=1$ if and only if $K(p,q)$ is a torus knot (or
link).
\end{rem}

\myparagraph{Main statements} The following will be established
below:

\begin{thm}\label{mainthm}
Let $K(p,q)$ be a given two-bridge knot (or link) and let
$(M_n(p,q))_{n=2}^\infty$ be a sequence of meridian-cyclic $n$-fold
branched coverings of $S^3$, branched along $K(p,q)$. Then:
\begin{equation}\label{upper:estimate}
c(M_n(p,q))\leqslant n(p-1)\quad\forall n.
\end{equation}
If in addition $K(p,q)$ is hyperbolic then the following inequality
holds for $n\geqslant 7$ with $c=4$:
\begin{equation}\label{lower:estimate}
c(M_n(p,q)) > n\cdot\left( 1 - \frac{c \pi^2}{n^2}\right)^{3/2}\cdot
\max \{ 2, 2 \ell(p,q) - 2.6667\ldots\ \};
\end{equation}
moreover, if $K(p,q)$ is neither $K(5,2)$ nor $K(7,3)$,
then the inequality holds for $n \geqslant 6$
with $c=2\sqrt{2}$.
\end{thm}

\begin{rem}
\emph{Combining the inequalities~(\ref{upper:estimate})
and~(\ref{lower:estimate}), and letting $n$ tend to infinity,
one gets the qualitative result that the complexity of
$M_n(p,q)$ is asymptotically equal to $n$ up to a multiplicative
constant.}
\end{rem}

To state our next result, we recall that the T-invariant $T(G)$ of a
finitely presented group $G$ was defined in~\cite{Delzant} as the
minimal number $t$ such that $G$ admits a presentation with $t$
relations of length $3$ and an arbitrary number of relations of
length at most $2$. A presentation with this property will be called
\emph{triangular}.

\begin{prop}\label{T:prop}
For $n\geqslant 2$ let $M_n(p,q)$
be a meridian-cyclic $n$-fold branched coverings of
$S^3$, branched along a two-bridge knot (or link) $K(p,q)$. Then:
$$T(\pi_1(M_n(p,q)))\leqslant n(p-1).$$
\end{prop}

We note that some connections between the Matveev complexity of a
closed 3-manifold and the T-invariant of its fundamental group
were already discussed in~\cite{pepe}.

The proofs of inequalities~(\ref{upper:estimate})
and~(\ref{lower:estimate}) are completely independent of each other.
We will begin with the latter, which is established in Section~2,
and then prove the former (together with Proposition~\ref{T:prop},
which follows from the same argument) in Section~3. The special case
where $\ell(p,q)$ equals 2 is discussed in detail in Section~4.

\myparagraph{Acknowledgment} We are grateful to David Futer for very
useful information and discussion on the results recently
obtained in~\cite{FKP} and~\cite{GF}.

\newpage

\section{Hyperbolic volume and the twist number:\\ The lower
estimate} \noindent We begin by recalling that a manifold is
hyperbolic if it has a Riemannian metric of constant sectional
curvature $-1$. We will use in the sequel many facts from hyperbolic
geometry without explicit reference, see for
instance~\cite{Adams,lectures,bible}.

The two versions of inequality~(\ref{lower:estimate}) are readily
deduced by combining the following three propositions. Here and
always in the sequel $v_3=1.01494\ldots$ denotes the volume of the
regular ideal tetrahedron in hyperbolic $3$-space $\matH^3$, and
``vol(M)'' is the hyperbolic volume of a manifold $M$.

\begin{prop}\label{obvious:prop}
If $M$ is a closed orientable hyperbolic manifold then
$$\vol(M)<c(M)\cdot v_3.$$
\end{prop}

\begin{prop}\label{limit:prop}
If $K(p,q)$ is hyperbolic then $M_n(p,q)$, as defined in the
statement of Theorem~\ref{mainthm}, is hyperbolic for $n \geqslant
4$. Moreover  the following inequality
holds for $n\geqslant 7$ with $c=4$:
\begin{equation}\label{quot:vol:lim}
\vol(M_n(p,q)) \geqslant n\cdot \left( 1 - \frac{c\pi^2}{n^2}
\right)^{3/2} \cdot \vol(S^3\setminus K(p,q)),
\end{equation}
and, if $K(p,q)$ is neither $K(5,2)$ nor $K(7,3)$,
then the inequality holds for $n \geqslant 6$
with $c=2\sqrt{2}$.
\end{prop}

\begin{prop}\label{vol:est:prop}
If $K(p,q)$ is hyperbolic then
\begin{equation}\label{vol:est:eq}
\vol(S^3\setminus K(p,q))\geqslant
v_3\cdot\max\{2,2\ell(p,q)-2.6667\ldots\ \}.
\end{equation}
\end{prop}

We begin proofs by establishing the general
connection between complexity and the hyperbolic volume:

\medskip
\noindent \emph{Proof of Proposition~\ref{obvious:prop}.} Set $k=c(M)$. Being
hyperbolic, $M$ is prime and not one of the exceptional manifolds
$S^3$, $\matR\matP^3$, $S^2\times S^1$, or
$L(3,1)$, so there exists a realization of $M$ as a gluing of $k$
tetrahedra. If $\Delta$ denotes the abstract tetrahedron, this
realization induces continuous maps $\sigma_i:\Delta\to M$ for
$i=1,\ldots,k$ given by the restrictions to the various tetrahedra
of the projection from the disjoint union of the tetrahedra to $M$.
Note that each $\sigma_i$ is injective on the interior of $\Delta$
but maybe not on the boundary. Since the gluings used to pair the
faces of the tetrahedra in the construction of $M$ are simplicial,
it follows that $\sum_{i=1}^k \sigma_i$ is a singular $3$-cycle,
which of course represents the fundamental class $[M]\in
H_3(M;\matZ)$.

We consider now the universal covering $\matH^3\to M$. Since
$\Delta$ is simply connected, it is possible to lift $\sigma_i$ to
a map $\widetilde{\sigma}_i:\Delta\to\matH^3$. We then define the
simplicial map $\widetilde{\tau}_i:\Delta\to\matH^3$ which agrees
with $\widetilde{\sigma}_i$ on the vertices, where geodesic convex
combinations are used in $\matH^3$ to define the notion of
``simplicial''. We also denote by $\tau_i:\Delta\to M$ the
composition of $\widetilde{\tau}_i$ with the projection
$\matH^3\to M$. It is immediate to see that $\sum_{i=1}^k \tau_i$
is again a singular $3$-cycle in $M$. Using this and taking convex
combinations in $\matH^3$, one can actually check that the cycles
$\sum_{i=1}^k \sigma_i$ and $\sum_{i=1}^k \tau_i$ are homotopic
to each other.
Therefore, since the first cycle represents $[M]$, the latter also
does, which implies that $\bigcup_{i=1}^k\tau_i(\Delta)$ is equal
to $M$, otherwise $\sum_{i=1}^k \tau_i$ would be homotopic to a
map with 2-dimensional image.

Next we note that $\widetilde{\tau}_i(\Delta)$ is a compact
geodesic tetrahedron in $\matH^3$, so its volume is less than
$v_3$, see~\cite{lectures}. Moreover the volume of
$\tau_i(\Delta)$ is at most equal to the volume of
$\widetilde{\tau}_i(\Delta)$, because the projection $\matH^3\to M$
is a local isometry, and the volume of $M$ is at most the sum of
the volumes of the $\tau_i(\Delta)$'s, because we have shown above
that $M$ is covered by the $\tau_i(\Delta)$'s (perhaps
with some overlapping). This establishes the proposition.
\finedimo

\medskip
\noindent \emph{Proof of Proposition}~\ref{limit:prop}. This is
actually a direct application of Theorem~3.5 of~\cite{FKP}. We only
need to note that in~\cite{FKP} the result is stated for hyperbolic
(not necessarily two-bridge) knots (rather than links), but it is
easy to see that the proof (based on~\cite{Adams2002} and
Theorem~1.1 of~\cite{FKP}) works well also for hyperbolic two-bridge
links and their meridian-cyclic coverings. \finedimo

Before getting to the proof of Proposition~\ref{vol:est:prop} we
establish the characterization of $\ell(p,q)$ stated in the first
section:

\medskip
\noindent \emph{Proof of Proposition~\ref{strange:ell:prop}.}
Under assumption~(\ref{pq:hp}), we know that the relevant pairs $(p',q')$
are those with $p'$ equal to $p$ and $q'$ equal to either $p-q$ or
$r$ or $p-r$, where $1\leqslant r\leqslant p-1$ and
$q\cdot r\equiv 1$ (mod $p$).

We begin by noting that if we take
positive continued fraction expansions of $p/q$ and $p/(p-q)$
we find $1$ as the first coefficient in one case and a number greater than $1$
in the other case. Supposing first that $p/q=[1,a_2,a_3,\ldots,a_k]$
it is now easy to see that $p/(p-q)=[a_2+1,a_3,\ldots,a_k]$, so the minimized positive
expansion of $p/(p-q)$ has length $k-1$. The same argument with switched roles shows that
if the first coefficient $a_1$ of the minimized positive expansion of $p/q$ is larger than $1$
then the length of the expansion of $p/(p-q)$ is $k+1$. Therefore the minimal
length we can obtain using $q$ and $p-q$ is indeed $\ell(p,q)$.

Supposing $p/q=[a_1,\ldots,a_k]$, we next choose
$s$ with $1\leqslant s\leqslant p-1$
and $q\cdot s\equiv(-1)^{k-1}$ (mod $p$), and we note that
$\{s,p-s\}=\{r,p-r\}$. Now it is not difficult to see that
$p/s$ has a positive continued fraction expansion
$p/s=[a_k,a_{k-1},\ldots,a_2,a_1]$. Note that this may or not be
a minimized expansion, depending on whether $a_1$ is greater than 1
or equal to 1, but the length of the minimized version
is $\ell(p,q)$ anyway, thanks to
the definition we have given. By the above argument, since $a_k>1$,
the length of the minimized positive expansion of $p/(p-s)$ is
1 more than that of $p/s$, and the proposition is established.
\finedimo

\medskip
\noindent \emph{Proof of Proposition}~\ref{vol:est:prop}. This will
be based on results of Cao-Meyer\-hoff~\cite{CM} and
Gu\'{e}ritaud-Futer~\cite{GF}. Note that~(\ref{vol:est:eq}) is
equivalent to the two inequalities
\begin{eqnarray}
\vol(S^3\setminus K(p,q)) & \geqslant & 2 v_3 \label{CM:est:eq} \\
\vol(S^3\setminus K(p,q)) & \geqslant & v_3 \cdot (2 \ell(p,q) -
2.6667\ldots\ ) \label{GF:est:ell}.
\end{eqnarray}
Now, Cao and Meyerhoff have proved in~\cite{CM} that the
figure-eight knot complement (namely $S^3 \setminus K(5,2)$ in our
notation) and its sibling manifold (which can be described as the
$(5,1)$-Dehn surgery on the right-handed Whitehead link)
are the orientable cusped hyperbolic $3$-mani\-folds of minimal volume, and
they are the only such $3$-manifolds. Each has volume equal to $2
v_3 = 2.02988\ldots$ , which implies inequality~(\ref{CM:est:eq})
directly.

To establish~(\ref{GF:est:ell}) we need to recall some terminology introduced by Lackenby
in~\cite{L}. A \emph{twist} in a link diagram $D\subset\matR^2$ is
either a maximal collection of bigonal regions of $\matR^2\setminus
D$ arranged in a row, or a single crossing with no incident bigonal
regions. The \emph{twist number} $t(D)$ of $D$ is the total number
of twists in $D$. Moreover $D$ is called \emph{twist-reduced} if it
is alternating and whenever $\gamma\subset\matR^2$ is a simple
closed curve meeting $D$ transversely at two crossing only, one of
the two portions into which $\gamma$ separates $D$ is contained in a
twist. (This is not quite the definition in~\cite{L}, but it is
easily recognized to be equivalent to it for alternating diagrams.)

Lackenby proved in~\cite{L} that if $D$ is a prime
twist-reduced diagram of a hyperbolic link $L$ in $S^3$ then
\begin{equation*}
v_3  \cdot (t(D) -2) \leqslant  \vol (S^3 \setminus L) \leqslant  10
\cdot v_3 \cdot (t(D) - 1) ,
\end{equation*}
where $v_3$ is the volume of the regular ideal tetrahedron. These
estimates were improved for the case of hyperbolic two-bridge links
by Gu{\'e}ritaud and Futer~\cite{GF}. More exactly, if $D$ is a
reduced alternating diagram of a hyperbolic two-bridge link $L$,
then by~\cite[Theorem~B.3]{GF}
\begin{equation}\label{GF:estimates}
2 v_3  \cdot t(D) - 2.7066\ldots\  <  \vol (S^3 \setminus L) < 2 v_8
\cdot (t(D) - 1) ,
\end{equation}
where $v_8$ is the volume of the regular ideal octahedron.

Using the first inequality in~(\ref{GF:estimates}), the next result
implies~(\ref{GF:est:ell}), which completes the proof of
Proposition~\ref{vol:est:prop} and hence of
inequality~(\ref{lower:estimate}) in Theorem~\ref{mainthm}:

\begin{prop}\label{twist:prop}
The link $K(p,q)$ has a twist-reduced diagram with twist number
$\ell(p,q)$.
\end{prop}

\begin{proof} The required diagram $D$
is simply given by the so-called Conway normal form of $K(p,q)$
associated to the minimized positive continued fraction expansion
$[a_1,\ldots,a_k]$ of $p/q$. The definition of the Conway normal
form differs for even and odd $k$, and it is described in
Fig.~\ref{conway:form:gen:fig}.
\begin{figure}
\begin{center}
\unitlength=0.26mm
\begin{picture}(390,60)(0,0)
\thicklines

\qbezier(20,0)(4,10)(20,20) \qbezier(20,40)(4,50)(20,60)

\qbezier(20,40)(30,30)(40,20)
\qbezier(20,20)(22,22)(28,28)
\qbezier(32,32)(34,34)(40,40)

\put(50,30){\circle*{2}}
\put(55,30){\circle*{2}}
\put(60,30){\circle*{2}}
\put(57,37){\makebox(0,0){$\scriptstyle a_1$}}
\qbezier(20,0)(30,0)(90,0)

\qbezier(70,40)(80,30)(90,20)
\qbezier(70,20)(72,22)(78,28)
\qbezier(82,32)(84,34)(90,40)

\qbezier(90,0)(100,10)(109,20)
\qbezier(90,20)(92,18)(98,12)
\qbezier(102,8)(104,6)(110,0)

\put(120,10){\circle*{2}}
\put(125,10){\circle*{2}}
\put(130,10){\circle*{2}}
\put(127,17){\makebox(0,0){$\scriptstyle a_2$}}
\qbezier(90,40)(95,40)(160,40)

\qbezier(140,0)(150,10)(160,20)
\qbezier(140,20)(142,18)(148,12)
\qbezier(152,8)(154,6)(160,0)

\qbezier(160,40)(170,30)(180,20) \qbezier(160,20)(162,22)(168,28)
\qbezier(172,32)(174,34)(180,40)

\put(190,30){\circle*{2}}
\put(195,30){\circle*{2}}
\put(200,30){\circle*{2}}
\put(197,37){\makebox(0,0){$\scriptstyle a_3$}}
\qbezier(160,0)(190,0)(230,0)

\qbezier(210,40)(220,30)(230,20)
\qbezier(210,20)(212,22)(218,28)
\qbezier(222,32)(224,34)(230,40)

\qbezier(20,60)(30,60)(230,60)

\put(255,0){\circle*{2}} \put(260,0){\circle*{2}}
\put(265,0){\circle*{2}}

\put(255,20){\circle*{2}} \put(260,20){\circle*{2}}
\put(265,20){\circle*{2}}

\put(255,40){\circle*{2}} \put(260,40){\circle*{2}}
\put(265,40){\circle*{2}}

\put(255,60){\circle*{2}}
\put(260,60){\circle*{2}}
\put(265,60){\circle*{2}}

\qbezier(290,0)(300,10)(309,20) \qbezier(290,20)(292,18)(298,12)
\qbezier(302,8)(304,6)(310,0)

\put(320,10){\circle*{2}} \put(325,10){\circle*{2}}
\put(330,10){\circle*{2}} \put(327,17){\makebox(0,0){$\scriptstyle
a_k$}} \qbezier(290,40)(295,40)(360,40)
\qbezier(290,60)(295,60)(360,60)

\qbezier(340,0)(350,10)(360,20) \qbezier(340,20)(342,18)(348,12)
\qbezier(352,8)(354,6)(360,0)

\qbezier(360,20)(376,30)(360,40) \qbezier(360,0)(387,7)(390,30)
\qbezier(360,60)(387,53)(390,30)
\end{picture}
\begin{picture}(360,80)(0,0)
\thicklines \qbezier(20,0)(4,10)(20,20) \qbezier(20,40)(4,50)(20,60)

\qbezier(20,40)(30,30)(40,20) \qbezier(20,20)(22,22)(28,28)
\qbezier(32,32)(34,34)(40,40)

\put(50,30){\circle*{2}} \put(55,30){\circle*{2}}
\put(60,30){\circle*{2}} \put(57,37){\makebox(0,0){$\scriptstyle
a_1$}} \qbezier(20,0)(30,0)(90,0)

\qbezier(70,40)(80,30)(90,20) \qbezier(70,20)(72,22)(78,28)
\qbezier(82,32)(84,34)(90,40)

\qbezier(90,0)(100,10)(109,20) \qbezier(90,20)(92,18)(98,12)
\qbezier(102,8)(104,6)(110,0)

\put(120,10){\circle*{2}} \put(125,10){\circle*{2}}
\put(130,10){\circle*{2}} \put(127,17){\makebox(0,0){$\scriptstyle
a_2$}} \qbezier(90,40)(95,40)(160,40)

\qbezier(140,0)(150,10)(160,20) \qbezier(140,20)(142,18)(148,12)
\qbezier(152,8)(154,6)(160,0)

\qbezier(160,40)(170,30)(180,20) \qbezier(160,20)(162,22)(168,28)
\qbezier(172,32)(174,34)(180,40)

\put(190,30){\circle*{2}} \put(195,30){\circle*{2}}
\put(200,30){\circle*{2}} \put(197,37){\makebox(0,0){$\scriptstyle
a_3$}} \qbezier(160,0)(190,0)(230,0)

\qbezier(210,40)(220,30)(230,20) \qbezier(210,20)(212,22)(218,28)
\qbezier(222,32)(224,34)(230,40)

\qbezier(20,60)(30,60)(230,60)

\put(255,0){\circle*{2}} \put(260,0){\circle*{2}}
\put(265,0){\circle*{2}}

\put(255,20){\circle*{2}} \put(260,20){\circle*{2}}
\put(265,20){\circle*{2}}

\put(255,40){\circle*{2}} \put(260,40){\circle*{2}}
\put(265,40){\circle*{2}}

\put(255,60){\circle*{2}} \put(260,60){\circle*{2}}
\put(265,60){\circle*{2}}

\qbezier(290,40)(300,30)(310,20) \qbezier(290,20)(292,22)(298,28)
\qbezier(302,32)(304,34)(310,40)

\put(320,30){\circle*{2}} \put(325,30){\circle*{2}}
\put(330,30){\circle*{2}} \put(327,37){\makebox(0,0){$\scriptstyle
a_k$}} \qbezier(290,0)(320,0)(360,0)
\qbezier(290,60)(320,60)(360,60)

\qbezier(340,40)(350,30)(360,20) \qbezier(340,20)(342,22)(348,28)
\qbezier(352,32)(354,34)(360,40)

\qbezier(360,0)(374,10)(360,20) \qbezier(360,40)(374,50)(360,60)
\end{picture}

\end{center}
\caption{The Conway normal form of a two-bridge link. The number of half-twists
of the appropriate type in the $j$-th portion of the diagram
is given by the positive integer $a_j$. The upper picture refers to
the case of even $k$ and the lower picture to the case of odd $k$.}
\label{conway:form:gen:fig}
\end{figure}
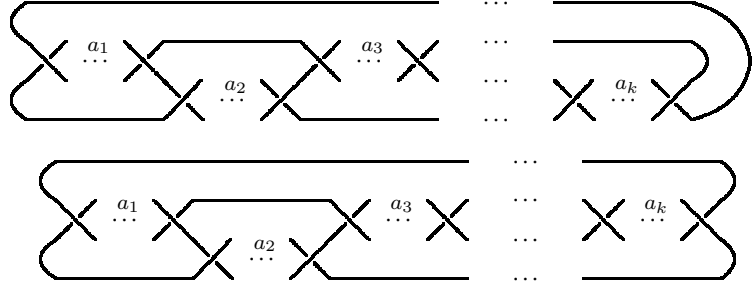
Two specific examples are also shown in
Fig.~\ref{conway:form:ex:fig}.
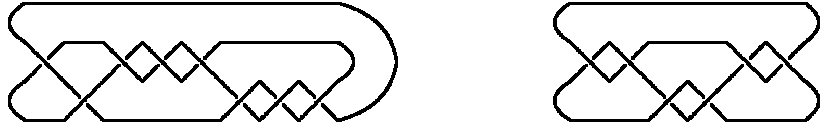
\begin{figure}
\begin{center}
\unitlength=0.26mm
\begin{picture}(210,60)(0,0)
\thicklines

\qbezier(20,0)(4,10)(20,20) \qbezier(20,40)(4,50)(20,60)
\qbezier(20,0)(30,0)(40,0)

\qbezier(20,40)(30,30)(40,20) \qbezier(20,20)(22,22)(28,28)
\qbezier(32,32)(34,34)(40,40)

\qbezier(40,0)(50,10)(60,20) \qbezier(40,20)(42,18)(48,12)
\qbezier(52,8)(54,6)(60,0) \qbezier(40,40)(50,40)(60,40)

\qbezier(60,40)(70,30)(80,20) \qbezier(60,20)(62,22)(68,28)
\qbezier(72,32)(74,34)(80,40)

\qbezier(80,40)(90,30)(100,20) \qbezier(80,20)(82,22)(88,28)
\qbezier(92,32)(94,34)(100,40)

\qbezier(100,40)(110,30)(120,20) \qbezier(100,20)(102,22)(108,28)
\qbezier(112,32)(114,34)(120,40)

\qbezier(60,0)(80,0)(120,0)

\qbezier(120,0)(130,10)(140,20) \qbezier(120,20)(122,18)(128,12)
\qbezier(132,8)(134,6)(140,0)

\qbezier(140,0)(150,10)(160,20) \qbezier(140,20)(142,18)(148,12)
\qbezier(152,8)(154,6)(160,0)

\qbezier(160,0)(170,10)(180,20) \qbezier(160,20)(162,18)(168,12)
\qbezier(172,8)(174,6)(180,0)

\qbezier(120,40)(150,40)(180,40) \qbezier(20,60)(50,60)(180,60)

\qbezier(180,20)(196,30)(180,40) \qbezier(180,0)(207,7)(210,30)
\qbezier(180,60)(207,53)(210,30)
\end{picture}
\qquad\qquad
\begin{picture}(160,60)(0,0)
\thicklines

\qbezier(20,0)(4,10)(20,20) \qbezier(20,40)(4,50)(20,60)

\qbezier(20,40)(30,30)(40,20) \qbezier(20,20)(22,22)(28,28)
\qbezier(32,32)(34,34)(40,40)

\qbezier(40,40)(50,30)(60,20) \qbezier(40,20)(42,22)(48,28)
\qbezier(52,32)(54,34)(60,40)

\qbezier(20,0)(40,0),(60,0)

\qbezier(60,0)(70,10)(80,20) \qbezier(60,20)(62,18)(68,12)
\qbezier(72,8)(74,6)(80,0)

\qbezier(80,0)(90,10)(100,20) \qbezier(80,20)(82,18)(88,12)
\qbezier(92,8)(94,6)(100,0)

\qbezier(100,40)(80,40),(60,40)

\qbezier(100,40)(110,30)(120,20) \qbezier(100,20)(102,22)(108,28)
\qbezier(112,32)(114,34)(120,40)

\qbezier(120,40)(130,30)(140,20) \qbezier(120,20)(122,22)(128,28)
\qbezier(132,32)(134,34)(140,40)

\qbezier(100,0)(120,0),(140,0) \qbezier(20,60)(30,60)(140,60)

\qbezier(140,0)(154,10)(140,20) \qbezier(140,40)(154,50)(140,60)
\end{picture}
\end{center}
\caption{Conway diagrams of $K(23,13)$ and $K(12,5)$. Note that the
required expansions are $23/13=[1,1,3,3]$ and $12/5=[2,2,2]$.}
\label{conway:form:ex:fig}
\end{figure}

Since the $a_j$'s are positive, it is quite obvious that the Conway
normal diagram $D$ always gives an alternating diagram, besides
being of course prime. The twists of this diagram are almost always
the obvious ones obtained by grouping together the first $a_1$
half-twists, then the next $a_2$, and so on. An exception has to be
made, however, when $a_1$ equals 1, because in this case the first
half-twist can be grouped with the next $a_2$ to give a single twist
(as in Fig.~\ref{conway:form:ex:fig}-left). Note that $a_k>1$ by
assumption, so no such phenomenon appears at the other end. Since
our definition of $\ell(p,q)$ is precisely $k$ if $a_1>1$ and $k-1$
if $a_1=1$, we see that indeed the diagram always has $\ell(p,q)$
twists.

Before proceeding we note that if $a_1=1$ then the Conway normal
form for $K(p,q)$ is actually the same, as a diagram, as the
mirror image of the Conway normal form for $K(p,p-q)$. The picture
showing this assertion gives a geometric proof of the fact that if
$p/q=[1,a_2,\ldots,a_k]$ then $p/(p-q)=[a_2+1,a_3,\ldots,a_k]$, used
in Proposition~\ref{strange:ell:prop}. So we can proceed assuming
that $a_1>1$. In particular, each bigonal region of $S^2\setminus D$
is one of the $(a_1-1)+(a_2-1)+\ldots+(a_k-1)$ created when
inserting the $a_1,a_2,\ldots,a_k$ half-twists of the normal form.

To prove that $D$ is twist-reduced, let us look for a curve
$\gamma$ as in the definition, namely one that intersects $D$
transversely at two crossings. Near each such intersection,
$\gamma$ must be either horizontal or vertical (see
Fig.~\ref{conway:form:gen:fig}). Let us first show that if it
meets some crossing $c$ of $D$ horizontally then $c$ is the
crossing arising from the single half-twist that corresponds to
some coefficient $a_j$ equal to 1. If this is not the case, then
either to the left or to the right of $c$ there is a bigonal
region of $S^2\setminus D$. Then $\gamma$ must meet horizontally
the crossing at the other end of this bigonal region, which
readily implies that $\gamma$ cannot meet the diagram in two
points only.

Having shown that $\gamma$ can only be vertical when it intersects
vertices, except at the vertices arising from the $a_j$'s with
$a_j=1$, let us give labels $R_0,R_1,\ldots,R_k,R_{k+1}$ to the
non-bigonal regions of $S^2\setminus D$, as in
Fig.~\ref{conway:reg:lab:fig},
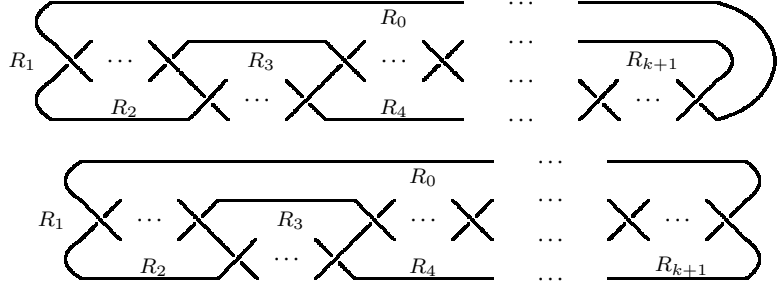
\begin{figure}
\begin{center}
\unitlength=0.26mm
\begin{picture}(390,60)(0,0)
\thicklines

\put(195,52){\makebox(0,0){$\scriptstyle R_0$}}
\put(5,30){\makebox(0,0){$\scriptstyle R_1$}}
\put(57,6){\makebox(0,0){$\scriptstyle R_2$}}
\put(127,30){\makebox(0,0){$\scriptstyle R_3$}}
\put(195,6){\makebox(0,0){$\scriptstyle R_4$}}
\put(327,30){\makebox(0,0){$\scriptstyle R_{k+1}$}}

\qbezier(20,0)(4,10)(20,20) \qbezier(20,40)(4,50)(20,60)

\qbezier(20,40)(30,30)(40,20) \qbezier(20,20)(22,22)(28,28)
\qbezier(32,32)(34,34)(40,40)

\put(50,30){\circle*{2}} \put(55,30){\circle*{2}}
\put(60,30){\circle*{2}} \qbezier(20,0)(30,0)(90,0)

\qbezier(70,40)(80,30)(90,20) \qbezier(70,20)(72,22)(78,28)
\qbezier(82,32)(84,34)(90,40)

\qbezier(90,0)(100,10)(109,20) \qbezier(90,20)(92,18)(98,12)
\qbezier(102,8)(104,6)(110,0)

\put(120,10){\circle*{2}} \put(125,10){\circle*{2}}
\put(130,10){\circle*{2}} \qbezier(90,40)(95,40)(160,40)

\qbezier(140,0)(150,10)(160,20) \qbezier(140,20)(142,18)(148,12)
\qbezier(152,8)(154,6)(160,0)

\qbezier(160,40)(170,30)(180,20) \qbezier(160,20)(162,22)(168,28)
\qbezier(172,32)(174,34)(180,40)

\put(190,30){\circle*{2}} \put(195,30){\circle*{2}}
\put(200,30){\circle*{2}} \qbezier(160,0)(190,0)(230,0)

\qbezier(210,40)(220,30)(230,20) \qbezier(210,20)(212,22)(218,28)
\qbezier(222,32)(224,34)(230,40)

\qbezier(20,60)(30,60)(230,60)

\put(255,0){\circle*{2}} \put(260,0){\circle*{2}}
\put(265,0){\circle*{2}}

\put(255,20){\circle*{2}} \put(260,20){\circle*{2}}
\put(265,20){\circle*{2}}

\put(255,40){\circle*{2}} \put(260,40){\circle*{2}}
\put(265,40){\circle*{2}}

\put(255,60){\circle*{2}} \put(260,60){\circle*{2}}
\put(265,60){\circle*{2}}

\qbezier(290,0)(300,10)(309,20) \qbezier(290,20)(292,18)(298,12)
\qbezier(302,8)(304,6)(310,0)

\put(320,10){\circle*{2}} \put(325,10){\circle*{2}}
\put(330,10){\circle*{2}} \qbezier(290,40)(295,40)(360,40)
\qbezier(290,60)(295,60)(360,60)

\qbezier(340,0)(350,10)(360,20) \qbezier(340,20)(342,18)(348,12)
\qbezier(352,8)(354,6)(360,0)

\qbezier(360,20)(376,30)(360,40) \qbezier(360,0)(387,7)(390,30)
\qbezier(360,60)(387,53)(390,30)
\end{picture}
\begin{picture}(360,80)(0,0)
\thicklines

\put(195,52){\makebox(0,0){$\scriptstyle R_0$}}
\put(5,30){\makebox(0,0){$\scriptstyle R_1$}}
\put(57,6){\makebox(0,0){$\scriptstyle R_2$}}
\put(127,30){\makebox(0,0){$\scriptstyle R_3$}}
\put(195,6){\makebox(0,0){$\scriptstyle R_4$}}
\put(327,6){\makebox(0,0){$\scriptstyle R_{k+1}$}}

\qbezier(20,0)(4,10)(20,20) \qbezier(20,40)(4,50)(20,60)

\qbezier(20,40)(30,30)(40,20) \qbezier(20,20)(22,22)(28,28)
\qbezier(32,32)(34,34)(40,40)

\put(50,30){\circle*{2}} \put(55,30){\circle*{2}}
\put(60,30){\circle*{2}} \qbezier(20,0)(30,0)(90,0)

\qbezier(70,40)(80,30)(90,20) \qbezier(70,20)(72,22)(78,28)
\qbezier(82,32)(84,34)(90,40)

\qbezier(90,0)(100,10)(109,20) \qbezier(90,20)(92,18)(98,12)
\qbezier(102,8)(104,6)(110,0)

\put(120,10){\circle*{2}} \put(125,10){\circle*{2}}
\put(130,10){\circle*{2}} \qbezier(90,40)(95,40)(160,40)

\qbezier(140,0)(150,10)(160,20) \qbezier(140,20)(142,18)(148,12)
\qbezier(152,8)(154,6)(160,0)

\qbezier(160,40)(170,30)(180,20) \qbezier(160,20)(162,22)(168,28)
\qbezier(172,32)(174,34)(180,40)

\put(190,30){\circle*{2}} \put(195,30){\circle*{2}}
\put(200,30){\circle*{2}} \qbezier(160,0)(190,0)(230,0)

\qbezier(210,40)(220,30)(230,20) \qbezier(210,20)(212,22)(218,28)
\qbezier(222,32)(224,34)(230,40)

\qbezier(20,60)(30,60)(230,60)

\put(255,0){\circle*{2}} \put(260,0){\circle*{2}}
\put(265,0){\circle*{2}}

\put(255,20){\circle*{2}} \put(260,20){\circle*{2}}
\put(265,20){\circle*{2}}

\put(255,40){\circle*{2}} \put(260,40){\circle*{2}}
\put(265,40){\circle*{2}}

\put(255,60){\circle*{2}} \put(260,60){\circle*{2}}
\put(265,60){\circle*{2}}

\qbezier(290,40)(300,30)(310,20) \qbezier(290,20)(292,22)(298,28)
\qbezier(302,32)(304,34)(310,40)

\put(320,30){\circle*{2}} \put(325,30){\circle*{2}}
\put(330,30){\circle*{2}} \qbezier(290,0)(320,0)(360,0)
\qbezier(290,60)(320,60)(360,60)

\qbezier(340,40)(350,30)(360,20) \qbezier(340,20)(342,22)(348,28)
\qbezier(352,32)(354,34)(360,40)

\qbezier(360,0)(374,10)(360,20) \qbezier(360,40)(374,50)(360,60)
\end{picture}
\end{center}
\caption{Labels for the non-bigonal regions of the complement of
a Conway normal diagram.}
\label{conway:reg:lab:fig}
\end{figure}
and let us construct a graph $\Gamma$ with vertices
$R_0,R_1,\ldots,R_k,R_{k+1}$ and an edge joining $R_i$ to $R_j$ for
each segment through a crossing of $D$ going from $R_i$ to $R_j$. By
assumption $\gamma$ must correspond to a length-2 cycle in $\Gamma$.
Now for odd $k$ the connections existing in $\Gamma$ are precisely
as follows:
\begin{itemize}
\item an $a_{2j-1}$-fold connection between $R_0$ and $R_{2j}$ for
$j=1,\ldots,(k+1)/2$;
\item an $a_{2j}$-fold connection between $R_1$ and $R_{2j+1}$ for $j=1,\ldots,(k-1)/2$;
\item a single connection between $R_j$ and $R_{j+2}$ if $2\leqslant j\leqslant k-1$ and $a_j=1$.
\end{itemize}

Then the only length-2 cycles are the evident ones
either between $R_0$ and some $R_{2j}$ or between $R_1$ and some $R_{2j+1}$,
and the curve $\gamma$ corresponding to one of these cycles does
bound a portion of a twist of $D$, as required by the definition of
twist-reduced diagram.

A similar analysis for even $k$ completes the proof.
\end{proof}

The proof of Proposition~\ref{vol:est:prop} is complete.\finedimo

\medskip
\noindent \emph{Proof of inequality}~(\ref{lower:estimate}).
Combining~(\ref{CM:est:eq}) and the first inequality
in~(\ref{GF:estimates}) with Proposition~\ref{twist:prop} we get
\begin{equation*}
v_3 \cdot \max\{2, 2 \ell(p,q) - 2.6667\ldots\ \} < \vol
(S^3 \setminus K(p,q)).
\end{equation*}
Together with~(\ref{quot:vol:lim}), this formula implies that
\begin{equation*}
\left(1-\frac{c\pi^2}{n^2} \right)^{3/2} \cdot v_3 \cdot \max\{2, 2
\ell(p,q) - 2.6667\ldots\ \} \cdot n < \vol(M_n(p,q))
\end{equation*}
with $c=4$ and $n\geqslant 7$ in general, and with $c=2\sqrt{2}$ and
$n\geqslant 6$ whenever $K(p,q)$ is neither $K(5,2)$ nor $K(7,3)$.
The conclusion now readily follows from
Proposition~\ref{obvious:prop}. \finedimo

\begin{rem} It was pointed out by Gu\'eritaud and Futer in~\cite{GF} that the lower
bound in~\emph{(\ref{GF:estimates})} is asymptotically sharp. But it is
numerically not very effective in some cases.  As an example we will
discuss below the case $p/q = k + 1/m$, where the lower bound given by
\emph{(\ref{GF:estimates})},
which translates into our~\emph{(\ref{GF:est:ell})}, is worse than the
Cao-Meyerhoff lower bound given
by~\emph{(\ref{CM:est:eq})}, since $\ell(p,q)=2$.
\end{rem}

\begin{rem}
On the basis of some computer experiments,
we conjecture that the Whitehead link
complement (namely $S^3 \setminus K(8,3)$, with $\vol(S^3 \setminus
K(8,3)) = 3.66386\ldots$ ) has the smallest vo\-lu\-me among all
two-bridge two-component links.
\end{rem}

\section{Minkus polyhedral schemes and triangulations:\\ The upper estimate}
\noindent The proof of~(\ref{upper:estimate}) and Proposition~\ref{T:prop}
is based on a realization of $M_n(p,q)$ as the quotient of a certain
polyhedron under a gluing of its faces. This construction extends
one that applies to lens spaces and it is originally due to
Minkus~\cite{Mi}. We will briefly review it here following~\cite{MuV}.

\begin{figure}
\begin{center}
\unitlength=0.26mm
\begin{picture}(400,300)(-200,10)
\thicklines

\qbezier(0,110)(0,110)(0,300) \qbezier(0,110)(0,110)(200,10)
\qbezier(0,110)(0,110)(200,210) \qbezier(0,110)(0,110)(-140,250)

\put(0,110){\circle*{6}} \put(0,130){\circle*{6}}
\put(0,170){\circle*{6}} \put(0,190){\circle*{6}}
\put(0,210){\circle*{6}} \put(0,250){\circle*{6}}
\put(0,270){\circle*{6}} \put(0,290){\circle*{6}}
\put(10,130){\makebox(0,0)[lc]{${}_1$}}
\put(10,170){\makebox(0,0)[lc]{${}_{p-q-1}$}}
\put(10,190){\makebox(0,0)[lc]{${}_{p-q}$}}
\put(10,210){\makebox(0,0)[lc]{${}_{p-q+1}$}}
\put(10,250){\makebox(0,0)[lc]{${}_{q-1}$}}
\put(10,270){\makebox(0,0)[lc]{${}_{q}$}}
\put(10,290){\makebox(0,0)[lc]{${}_{p-1}$}}

\put(20,100){\circle*{6}} \put(60,80){\circle*{6}}
\put(80,70){\circle*{6}} \put(100,60){\circle*{6}}
\put(140,40){\circle*{6}} \put(160,30){\circle*{6}}
\put(180,20){\circle*{6}}

\put(10,100){\makebox(0,0)[rc]{${}_1$}}
\put(50,80){\makebox(0,0)[rc]{${}_{p-q-1}$}}
\put(70,70){\makebox(0,0)[rc]{${}_{p-q}$}}
\put(90,60){\makebox(0,0)[rc]{${}_{p-q+1}$}}
\put(130,40){\makebox(0,0)[rc]{${}_{q-1}$}}
\put(150,30){\makebox(0,0)[rc]{${}_{q}$}}
\put(170,20){\makebox(0,0)[rc]{${}_{p-1}$}}

\put(20,120){\circle*{6}} \put(60,140){\circle*{6}}
\put(80,150){\circle*{6}} \put(100,160){\circle*{6}}
\put(140,180){\circle*{6}} \put(160,190){\circle*{6}}
\put(180,200){\circle*{6}}

\put(20,110){\makebox(0,0)[cc]{${}_1$}}
\put(60,130){\makebox(0,0)[cc]{${}_{p-q-1}$}}
\put(80,140){\makebox(0,0)[cc]{${}_{p-q}$}}
\put(105,150){\makebox(0,0)[cc]{${}_{p-q+1}$}}
\put(140,170){\makebox(0,0)[cc]{${}_{q-1}$}}
\put(160,180){\makebox(0,0)[cc]{${}_{q}$}}
\put(180,190){\makebox(0,0)[cc]{${}_{p-1}$}}

\put(-20,130){\circle*{6}} \put(-50,160){\circle*{6}}
\put(-60,170){\circle*{6}} \put(-70,180){\circle*{6}}
\put(-90,200){\circle*{6}} \put(-100,210){\circle*{6}}
\put(-130,240){\circle*{6}}

\put(-30,130){\makebox(0,0)[rc]{${}_1$}}
\put(-60,160){\makebox(0,0)[rc]{${}_{p-q-1}$}}
\put(-70,170){\makebox(0,0)[rc]{${}_{p-q}$}}
\put(-80,180){\makebox(0,0)[rc]{${}_{p-q+1}$}}
\put(-100,200){\makebox(0,0)[rc]{${}_{q-1}$}}
\put(-110,210){\makebox(0,0)[rc]{${}_{q}$}}
\put(-140,240){\makebox(0,0)[rc]{${}_{p-1}$}}

\qbezier(80,70)(200,100)(160,190) \qbezier(80,150)(100,250)(0,270)
\qbezier(0,190)(-40,250)(-100,210)
\put(60,210){\makebox(0,0)[cc]{${R_i}$}}
\put(100,250){\makebox(0,0)[cc]{${R_{i+1}^{\prime}}$}}
\put(130,120){\makebox(0,0)[cc]{${R_{i+1}}$}}
\put(190,120){\makebox(0,0)[cc]{${R_{i+2}^{\prime}}$}}
\put(-40,200){\makebox(0,0)[cc]{${R_{i-1}}$}}
\put(-50,250){\makebox(0,0)[cc]{${R_{i}^{\prime}}$}}
\put(10,50){\circle*{6}} \put(-50,70){\circle*{6}}
\put(-70,120){\circle*{6}}

\put(-15,105){\makebox(0,0)[cc]{$N$}}
\put(160,205){\makebox(0,0)[cc]{$P_{i+1}$}}
\put(-15,270){\makebox(0,0)[cc]{$P_{i}$}}
\put(-95,230){\makebox(0,0)[cc]{$P_{i-1}$}}

\put(-15,185){\makebox(0,0)[cc]{$P'_{i}$}}
\put(65,165){\makebox(0,0)[cc]{$P'_{i+1}$}}
\put(85,85){\makebox(0,0)[cc]{$P'_{i+2}$}}
\end{picture}
\end{center} \caption{The Minkus polyhedral scheme for $M_n(p,q)$.}
\label{Minkus:fig}
\end{figure}
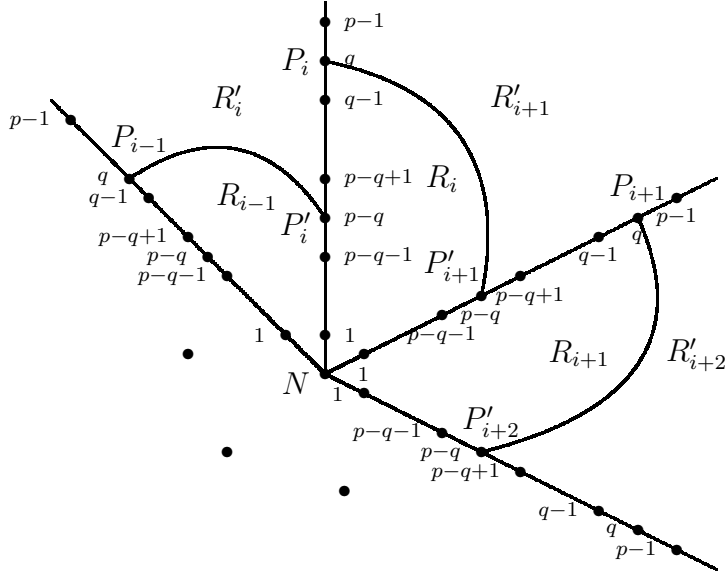

Let us begin from the case where $K(p,q)$ is a knot, \emph{i.e.}~$p$
is odd, whence $M_n(p,q)$ is uniquely defined by $p,q,n$. Recall
that by the assumption~(\ref{pq:hp}), $0<q<p$. Then we consider the
3-ball
$$B^3 = \{ (x,y,z) \in \matR^3:\ x^2 + y^2 +z^2 \leqslant 1 \}$$
and we draw on its boundary $n$ equally spaced great
semicircles joining the North pole $N=(0,0,1)$ to the South pole
$S=(0,0,-1)$. This decomposes $\partial B^3$ into $n$ cyclically
arranged congruent lunes $L_0,\ldots,L_{n-1}$. Now we insert $p-1$
equally spaced vertices on each semicircle, thus subdividing it into
$p$ identical segments, which allows us to view each lune $L_i$ as a
curvilinear polygon with $2p$ edges. Next, we denote by $P_i$ the
vertex on the semicircle $L_i\cap L_{i-1}$ which is $q$ segments
down from $N$, and by $P'_i$ the vertex which is $q$ segments up
from $S$ (indices are always meant modulo $n$). We then draw inside
$L_i$ an arc of great semicircle joining $P_i$ to $P'_{i+1}$, thus
bisecting $L_i$ into two regions that we denote by $R_i$ and
$R'_{i+1}$, with $R_i$ incident to $N$ and $R'_{i+1}$ incident to
$S$. Fig.~\ref{Minkus:fig} illustrates
the resulting decomposition of $\partial B^3$, which is
represented as $\matR^2 \cup \{\infty\}$ with $S=\infty$. In the
picture we assume $q>p/2$.

Summing up, we have subdivided $\partial B^3$ into $2n$
curvilinear polygons $R_i,R'_i$ for $i=0,\ldots, n-1$, each having
$p+1$ edges. The polygons $R_i$ are around $N$ and the polygons
$R'_i$ are around $S$, and there is a marked vertex $P_i$ shared by
$R_i$ and $R'_{i+1}$ (we will not need to use $P'_i$ again). It is
now possible to show that the manifold $M_n(p,q)$ is obtained from
$B^3$ by identifying $R_i$ with $R'_i$ on $\partial B^3$ for
$i=0, \ldots , n-1$ through an orientation-reversing simplicial
homeomorphism which matches the vertex $P_i$ of $R_i$ with the
vertex $P_{i-1}$ of $R'_i$.

\begin{figure}
\begin{center}
\unitlength=0.3mm
\begin{picture}(200,160)(-100,0)
\thicklines

\qbezier(0,60)(0,60)(0,160) \qbezier(0,60)(0,60)(100,10)
\qbezier(0,60)(0,60)(-100,10)

\put(0,60){\circle*{4}} \put(0,80){\circle*{4}}
\put(0,100){\circle*{4}} \put(0,120){\circle*{4}}
\put(0,140){\circle*{4}}

\put(10,80){\makebox(0,0)[cc]{${}_1$}}
\put(10,100){\makebox(0,0)[cc]{${}_2$}}
\put(10,125){\makebox(0,0)[cc]{${}_3$}}
\put(10,140){\makebox(0,0)[cc]{${}_4$}}
\put(-10,120){\makebox(0,0)[cc]{$P_1$}}

\put(20,50){\circle*{4}} \put(40,40){\circle*{4}}
\put(60,30){\circle*{4}} \put(80,20){\circle*{4}}

\put(20,40){\makebox(0,0)[cc]{${}_1$}}
\put(40,30){\makebox(0,0)[cc]{${}_2$}}
\put(60,20){\makebox(0,0)[cc]{${}_3$}}
\put(80,10){\makebox(0,0)[cc]{${}_4$}}
\put(60,40){\makebox(0,0)[cc]{$P_2$}}

\put(-20,50){\circle*{4}} \put(-40,40){\circle*{4}}
\put(-60,30){\circle*{4}} \put(-80,20){\circle*{4}}

\put(-20,60){\makebox(0,0)[cc]{${}_1$}}
\put(-40,50){\makebox(0,0)[cc]{${}_2$}}
\put(-63,38){\makebox(0,0)[cc]{${}_3$}}
\put(-80,30){\makebox(0,0)[cc]{${}_4$}}
\put(-60,20){\makebox(0,0)[cc]{$P_3$}}

\qbezier(40,40)(50,100)(0,120) \qbezier(60,30)(10,0)(-40,40)
\qbezier(0,100)(-50,100)(-60,30)

\put(25,75){\makebox(0,0)[cc]{${R_1}$}}
\put(60,100){\makebox(0,0)[cc]{${R_{2}^{\prime}}$}}
\put(0,40){\makebox(0,0)[cc]{${R_2}$}}
\put(-10,5){\makebox(0,0)[cc]{${R_{3}^{\prime}}$}}
\put(-25,75){\makebox(0,0)[cc]{${R_3}$}}
\put(-70,90){\makebox(0,0)[cc]{${R_{1}^{\prime}}$}}
\end{picture}
\end{center} \caption{The Minkus polyhedral scheme for $M_3(5,3)$.} \label{HW:fig}
\end{figure}
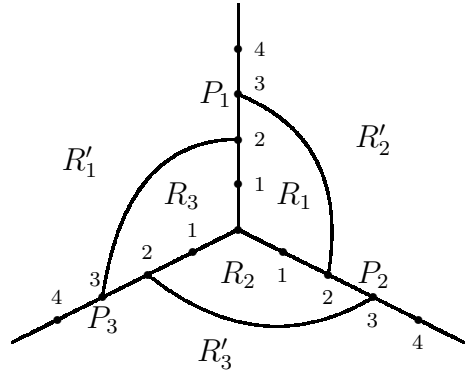

As an example, Fig.~\ref{HW:fig} shows the Minkus polyhedral
construction of the Hantzsche-Wendt manifold, that is $M_3(5,3)$ in
our notation.

\medskip\noindent
\emph{Proof of inequality}~(\ref{upper:estimate}) \emph{for odd $p$}. Referring
to the above polyhedral construction of $M_n(p,q)$, we subdivide
each $R_i$ into $p-1$ triangles by taking diagonals from the North
pole $N$, and each $R'_i$ so that the gluing between $R_i$ and
$R'_i$ matches the subdivision. Note that the ``diagonals'' are only
meant in a combinatorial sense, they cannot be taken as arcs of
great circles. We can then take (combinatorial) cones with vertex at
$N$ and bases at the triangles contained in the $R'_i$, thus getting
a subdivision of $B^3$ into $n(p-1)$ tetrahedra. By construction
the gluings on $\partial B^3$ restrict to gluings of the faces of
these tetrahedra, therefore $M_n(p,q)$ has a (loose) triangulation
made of $n(p-1)$ tetrahedra, and the proof is now complete. \finedimo

\begin{figure}
\begin{center}
\unitlength=0.26mm
\begin{picture}(400,300)(-200,10)
\thicklines

\qbezier(0,110)(0,110)(0,300) \qbezier(0,110)(0,110)(200,10)
\qbezier(0,110)(0,110)(200,210) \qbezier(0,110)(0,110)(-140,250)

\put(0,110){\circle*{6}} \put(0,130){\circle*{6}}
\put(0,170){\circle*{6}} \put(0,190){\circle*{6}}
\put(0,210){\circle*{6}} \put(0,250){\circle*{6}}
\put(0,270){\circle*{6}} \put(0,290){\circle*{6}}
\put(10,130){\makebox(0,0)[lc]{${}_1$}}
\put(10,170){\makebox(0,0)[lc]{${}_{p-q-1}$}}
\put(10,190){\makebox(0,0)[lc]{${}_{p-q}$}}
\put(10,210){\makebox(0,0)[lc]{${}_{p-q+1}$}}
\put(10,250){\makebox(0,0)[lc]{${}_{q-1}$}}
\put(10,270){\makebox(0,0)[lc]{${}_{q}$}}
\put(10,290){\makebox(0,0)[lc]{${}_{p-1}$}}

\put(20,100){\circle*{6}} \put(60,80){\circle*{6}}
\put(80,70){\circle*{6}} \put(100,60){\circle*{6}}
\put(140,40){\circle*{6}} \put(160,30){\circle*{6}}
\put(180,20){\circle*{6}}

\put(10,100){\makebox(0,0)[rc]{${}_1$}}
\put(50,80){\makebox(0,0)[rc]{${}_{p-q-1}$}}
\put(70,70){\makebox(0,0)[rc]{${}_{p-q}$}}
\put(90,60){\makebox(0,0)[rc]{${}_{p-q+1}$}}
\put(130,40){\makebox(0,0)[rc]{${}_{q-1}$}}
\put(150,30){\makebox(0,0)[rc]{${}_{q}$}}
\put(170,20){\makebox(0,0)[rc]{${}_{p-1}$}}

\put(20,120){\circle*{6}} \put(60,140){\circle*{6}}
\put(80,150){\circle*{6}} \put(100,160){\circle*{6}}
\put(140,180){\circle*{6}} \put(160,190){\circle*{6}}
\put(180,200){\circle*{6}}

\put(20,110){\makebox(0,0)[cc]{${}_1$}}
\put(60,130){\makebox(0,0)[cc]{${}_{p-q-1}$}}
\put(80,140){\makebox(0,0)[cc]{${}_{p-q}$}}
\put(105,150){\makebox(0,0)[cc]{${}_{p-q+1}$}}
\put(140,170){\makebox(0,0)[cc]{${}_{q-1}$}}
\put(160,180){\makebox(0,0)[cc]{${}_{q}$}}
\put(180,190){\makebox(0,0)[cc]{${}_{p-1}$}}

\put(-20,130){\circle*{6}} \put(-50,160){\circle*{6}}
\put(-60,170){\circle*{6}} \put(-70,180){\circle*{6}}
\put(-90,200){\circle*{6}} \put(-100,210){\circle*{6}}
\put(-120,230){\circle*{6}}

\put(-30,130){\makebox(0,0)[rc]{${}_1$}}
\put(-60,160){\makebox(0,0)[rc]{${}_{p-q-1}$}}
\put(-70,170){\makebox(0,0)[rc]{${}_{p-q}$}}
\put(-80,180){\makebox(0,0)[rc]{${}_{p-q+1}$}}
\put(-100,200){\makebox(0,0)[rc]{${}_{q-1}$}}
\put(-110,210){\makebox(0,0)[rc]{${}_{q}$}}
\put(-130,230){\makebox(0,0)[rc]{${}_{p-1}$}}

\qbezier(80,70)(180,120)(160,190)

\qbezier(80,150)(80,250)(0,270)

\qbezier(0,190)(-50,250)(-100,210)

\put(50,190){\makebox(0,0)[cc]{${R_i}$}}
\put(100,250){\makebox(0,0)[cc]{${R_{i+m}^{\prime}}$}}
\put(100,120){\makebox(0,0)[cc]{${R_{i+1}}$}}
\put(190,120){\makebox(0,0)[cc]{${R_{i+m+1}^{\prime}}$}}
\put(-30,190){\makebox(0,0)[cc]{${R_{i-1}}$}}
\put(-50,250){\makebox(0,0)[cc]{${R_{i+m-1}^{\prime}}$}}
\put(10,50){\circle*{6}} \put(-50,70){\circle*{6}}
\put(-70,120){\circle*{6}}
\end{picture}
\end{center} \caption{The Minkus polyhedral scheme for $M_{n,m}(p,q)$.}
\label{Minkus:k:fig}
\end{figure}
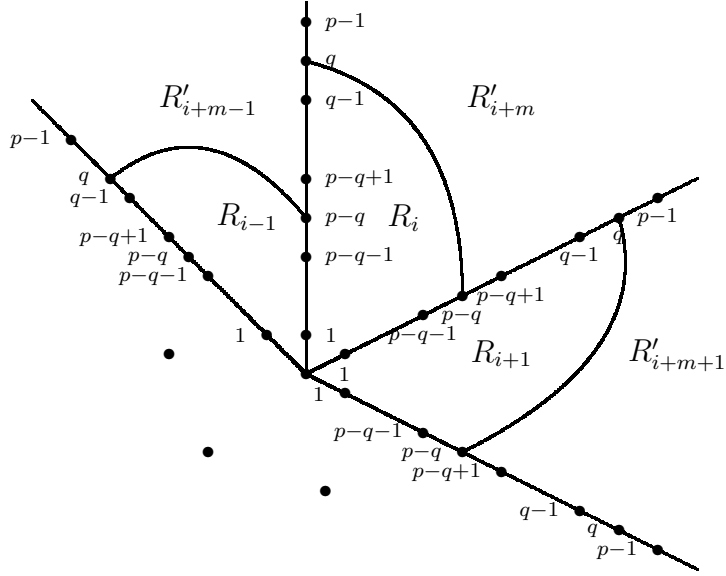

\medskip\noindent
\emph{Proof of inequality}~(\ref{upper:estimate}) \emph{for even $p$}. To
establish~(\ref{upper:estimate}) for even $p$, \emph{i.e.}~for
2-component two-bridge links, we extend to this case the Minkus
polyhedral construction, see~\cite{MuV}. The way to do this is
actually straight-forward: to realize the meridian-cyclic covering
$M_{n,m}(p,q)$ of $S^3$ bran\-ched along $K(p,q)$ we subdivide
$\partial B^3$ precisely as above, but we denote by $R_i$ and
$R'_{i+m}$ the two regions into which the lune $L_i$ is bisected.
Then we glue  $R_i$ to $R'_i$ by an orientation-reversing simplicial
homeomorphism matching the vertex $P_i$ of $R_i$ with the vertex
$P_{i-m}$ of $R'_i$. This construction is illustrated in
Fig.~\ref{Minkus:k:fig}. This realization of $M_{n,m}(p,q)$ again
induces a triangulation with $n(p-1)$ tetrahedra, which
proves~(\ref{upper:estimate}) also in this case. \finedimo

\medskip
\noindent \emph{Proof of Proposition}~\ref{T:prop}. Let us carry
out only the ``first half'' of the subdivision we did above of the
Minkus polyhedral realization of $M_n(p,q)$. Namely, we subdivide
the regions $R_i,R'_i$ on $\partial B^3$ into triangles, but then
we do not add anything inside $B^3$. This yields a
cellularization of $M_n(p,q)$ with 2-cells being triangles and with
a single $3$-cell. Therefore there is a triangular presentation of
$\pi_1(M_n(p,q))$ with precisely the same number of relations as the
number of triangles in this cellularization. And this number is
$n(p-1)$, because there are $2n(p-1)$ triangles on
$\partial B^3$, but they get glued in pairs. \finedimo

\section{An Example: $p/q = k + 1/m$}

\noindent As already noticed, the lower bound on the volume given
by~(\ref{GF:est:ell}) does not seem to provide very effective
estimates in some instances.  We consider in this section the case
where $p/q$ has the form $k + 1/m=(km+1)/m$ and $K(p,q)$ is
hyperbolic. Then $\ell(p,q)=2$, so the estimate provided
by~(\ref{GF:est:ell}) is weaker than that in~(\ref{CM:est:eq}).
However, denoting by $\mathcal B$ the Borromean rings, a hyperbolic
3-component link (depicted below) with volume $7.32772\ldots$ , we
can establish the following:

\begin{prop}\label{l2:prop}
If $k,m>1$ then $K(km+1,m)$ is hyperbolic and
$$\lim_{k,m\to\infty} \vol (S^3 \setminus K(km+1,m)) = \vol(S^3 \setminus \mathcal B) =
7.32772\ldots\ .$$
\end{prop}

\begin{proof}
We set $p=km+1$ and $q=m$ and note that $K(p,q)$ is hyperbolic by the
discussion in Section~1. For the limit of the volume,
the proof is always based on hyperbolic surgery, but it differs according
to whether $k$ and $m$ are even or odd.

\smallskip \noindent\textsc{Case 1: $k=2i$ and $m=2j$.} In this
case $K(p,q) = K(4ij+1, 2j)$, shown in
Fig.~\ref{Borromean:fig}-left, can be obtained by Dehn surgery on
the Borromean rings $\mathcal B = 6^3_2$ pictured in
Fig.~\ref{Borromean:fig}-right, doing a $(1/i)$-surgery on the first
small circle and a $(1/j)$-surgery on the another one,
see~\cite[Ch.~9]{Rolfsen} for details. Labels for links, here and in
the sequel, are also taken from~\cite{Rolfsen}. Thurston's
hyperbolic Dehn surgery theorem then implies that
$$\lim_{i,j\to\infty}
\vol(S^3 \setminus K(4ij+1,2j)) = \vol(S^3 \setminus \mathcal B) =
7.32772\ldots$$
as required.

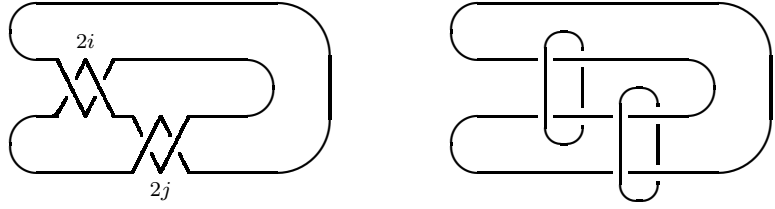
\begin{figure}
\begin{center}
\unitlength=.25mm
\begin{picture}(0,110)(0,30)
\put(-120,0){
\begin{picture}(0,130)(-40,0)
\thicklines

\qbezier(-100,40)(-100,40)(-55,40) \qbezier(-25,40)(-25,40)(0,40)
\qbezier(-100,70)(-100,70)(-95,70) \qbezier(-65,70)(-65,70)(-55,70)
\qbezier(-25,70)(-25,70)(0,70) \qbezier(-100,100)(-100,100)(-95,100)
\qbezier(-65,100)(-65,100)(0,100)
\qbezier(-100,130)(-100,130)(0,130)

\put(0,85){\oval(40,30)[r]} \put(0,85){\oval(100,90)[r]}
\put(-100,55){\oval(40,30)[l]} \put(-100,115){\oval(40,30)[l]}

\qbezier(-95,100)(-95,100)(-80,70) \qbezier(-97,70)(-95,70)(-90,80)
\qbezier(-80,100)(-80,100)(-85,90)
\qbezier(-80,100)(-80,100)(-65,70) \qbezier(-80,70)(-80,70)(-75,80)
\qbezier(-65,100)(-65,100)(-70,90) \qbezier(-55,40)(-55,40)(-40,70)
\qbezier(-55,70)(-55,70)(-50,60) \qbezier(-40,40)(-40,40)(-45,50)
\qbezier(-40,40)(-40,40)(-25,70) \qbezier(-25,40)(-25,40)(-30,50)
\qbezier(-40,70)(-40,70)(-35,60)

\put(-80,110){\makebox(0,0)[cc]{$\scriptstyle 2i$}}
\put(-40,30){\makebox(0,0)[cc]{$\scriptstyle 2j$}}
\end{picture}}
\put(120,0){\begin{picture}(0,130)(-40,0) \thicklines

\qbezier(-100,40)(-100,40)(-55,40) \qbezier(-25,40)(-25,40)(0,40)
\qbezier(-100,70)(-100,70)(-95,70) \qbezier(-65,70)(-65,70)(-55,70)
\qbezier(-25,70)(-25,70)(0,70) \qbezier(-100,100)(-100,100)(-95,100)
\qbezier(-65,100)(-65,100)(0,100)
\qbezier(-100,130)(-100,130)(0,130)

\put(0,85){\oval(40,30)[r]} \put(0,85){\oval(100,90)[r]}
\put(-100,55){\oval(40,30)[l]} \put(-100,115){\oval(40,30)[l]}

\put(-20,35){\oval(20,20)[b]} \put(-20,75){\oval(20,20)[t]}
\qbezier(-30,35)(-30,35)(-30,75) \qbezier(-10,45)(-10,45)(-10,65)

\put(-60,65){\oval(20,20)[b]} \put(-60,105){\oval(20,20)[t]}
\qbezier(-70,65)(-70,65)(-70,105) \qbezier(-50,75)(-50,75)(-50,95)

\qbezier(-75,40)(-75,40)(-35,40) \qbezier(-25,40)(-25,40)(0,40)
\qbezier(-95,70)(-95,70)(-75,70) \qbezier(-65,70)(-65,70)(-35,70)
\qbezier(-95,100)(-95,100)(-75,100)
\end{picture}
}
\end{picture} \caption{$K(4ij+1,2j)$
and $\mathcal B=6^3_2$.} \label{Borromean:fig}
\end{center}
\end{figure}

\smallskip \noindent\textsc{Case 2: $k=2i+1$ and $m=2j$.} In
this case $K(p,q) = K(4ij+2j+1, 2j)$, shown in
Fig.~\ref{Borromean:prime:fig}-left, can be obtained by Dehn surgery
on the $3$-component link $\mathcal B^{\prime} = 8^3_9$ pictured in
Fig.~\ref{Borromean:prime:fig}-right, doing a $(1/i)$-surgery on the
first small circle and a $(1/j)$-surgery the another one. By the
same arguments as in~\cite[p.~269-270]{lectures} we have  $\vol(S^3
\setminus \mathcal B^{\prime}) = \vol(S^3 \setminus \mathcal B)$.
Hence again
$$\lim_{i,j\to\infty} \vol(S^3 \setminus K(4ij+2j+1,2j)) = \vol(S^3
\setminus \mathcal B) = 7.32772\ldots\ .$$

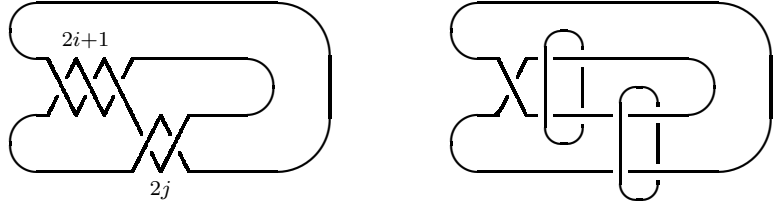
\begin{figure}[ht]
\begin{center}
\unitlength=.25mm
\begin{picture}(0,110)(0,30)
\put(-120,0){
\begin{picture}(0,130)(-40,0)
\thicklines

\qbezier(-100,40)(-100,40)(-55,40) \qbezier(-25,40)(-25,40)(0,40)
\qbezier(-25,70)(-25,70)(0,70) \qbezier(-55,100)(-55,100)(0,100)
\qbezier(-100,130)(-100,130)(0,130) \put(0,85){\oval(40,30)[r]}
\put(0,85){\oval(100,90)[r]} \put(-100,55){\oval(40,30)[l]}
\put(-100,115){\oval(40,30)[l]}

\qbezier(-100,100)(-100,100)(-85,70)
\qbezier(-100,70)(-100,70)(-95,80)
\qbezier(-85,100)(-85,100)(-90,90)
\qbezier(-85,100)(-85,100)(-70,70) \qbezier(-85,70)(-85,70)(-80,80)
\qbezier(-70,100)(-70,100)(-75,90)
\qbezier(-70,100)(-70,100)(-55,70) \qbezier(-70,70)(-70,70)(-65,80)
\qbezier(-55,100)(-55,100)(-60,90) \qbezier(-55,40)(-55,40)(-40,70)
\qbezier(-55,70)(-55,70)(-50,60) \qbezier(-40,40)(-40,40)(-45,50)
\qbezier(-40,40)(-40,40)(-25,70) \qbezier(-25,40)(-25,40)(-30,50)
\qbezier(-40,70)(-40,70)(-35,60)

\put(-80,110){\makebox(0,0)[cc]{$\scriptstyle 2i+1$}}
\put(-40,30){\makebox(0,0)[cc]{$\scriptstyle 2j$}}
\end{picture}}
\put(120,0){\begin{picture}(0,130)(-40,0) \thicklines

\qbezier(-100,40)(-100,40)(-55,40) \qbezier(-25,40)(-25,40)(0,40)
\qbezier(-100,70)(-100,70)(-95,70) \qbezier(-65,70)(-65,70)(-55,70)
\qbezier(-25,70)(-25,70)(0,70) \qbezier(-100,100)(-100,100)(-95,100)
\qbezier(-65,100)(-65,100)(0,100)
\qbezier(-100,130)(-100,130)(0,130)

\put(0,85){\oval(40,30)[r]} \put(0,85){\oval(100,90)[r]}
\put(-100,55){\oval(40,30)[l]} \put(-100,115){\oval(40,30)[l]}

\put(-20,35){\oval(20,20)[b]} \put(-20,75){\oval(20,20)[t]}
\qbezier(-30,35)(-30,35)(-30,75) \qbezier(-10,45)(-10,45)(-10,65)

\put(-60,65){\oval(20,20)[b]} \put(-60,105){\oval(20,20)[t]}
\qbezier(-70,65)(-70,65)(-70,105) \qbezier(-50,75)(-50,75)(-50,95)

\qbezier(-95,100)(-95,100)(-80,70) \qbezier(-97,70)(-95,70)(-90,80)
\qbezier(-80,100)(-80,100)(-85,90)

\qbezier(-75,40)(-75,40)(-35,40) \qbezier(-25,40)(-25,40)(0,40)
\qbezier(-80,70)(-80,70)(-75,70) \qbezier(-65,70)(-65,70)(-35,70)
\qbezier(-80,100)(-80,100)(-75,100)
\end{picture}}
\end{picture} \caption{$K(4ij+2j+1,2j)$ and $\mathcal B^{\prime} = 8^3_9$.}
\label{Borromean:prime:fig}
\end{center}
\end{figure}

\smallskip \noindent\textsc{Case 3: $k=2i$ and $m=2j+1$.} In this
case $K(p,q) = K(4ij+2i+1, 2j+1)$ is equivalent to
$K(p,q')=K(4ij+2i+1,2i)$ because $(2j+1)(2i) \equiv -1$ (mod $p$).
The conclusion then follows from Case 2.

\smallskip \noindent\textsc{Case 4: $k=2i+1$ and $m=2j+1$.} In this
case $K(p,q) = K(4ij+2i+2j+2, 2j+1)$, shown in
Fig.~\ref{Borromean:2prime:fig}-left, is a $2$-component link that
can be obtained by Dehn surgery on the two small circles of the
$4$-component link $\mathcal B^{\prime \prime}$ pictured in
Fig.~\ref{Borromean:2prime:fig}-right, doing a $(1/i)$-surgery and a
$(1/j)$-surgery. Since as above $\vol(S^3 \setminus \mathcal
B^{\prime \prime}) = \vol(S^3 \setminus \mathcal B)$, we have
$$\lim_{i,j\to\infty} \vol(S^3 \setminus K(4ij+2i+2j+2,2j+1)) =
7.32772\ldots$$ again.

\begin{figure}[ht]
\begin{center}
\unitlength=.25mm
\begin{picture}(0,110)(0,30)
\put(-120,0){
\begin{picture}(0,130)(-40,0)
\thicklines

\qbezier(-100,40)(-100,40)(-55,40) \qbezier(-10,40)(-10,40)(0,40)
\qbezier(-10,70)(-10,70)(0,70) \qbezier(-55,100)(-55,100)(0,100)
\qbezier(-100,130)(-100,130)(0,130) \put(0,85){\oval(40,30)[r]}
\put(0,85){\oval(100,90)[r]} \put(-100,55){\oval(40,30)[l]}
\put(-100,115){\oval(40,30)[l]}

\qbezier(-100,100)(-100,100)(-85,70)
\qbezier(-100,70)(-100,70)(-95,80)
\qbezier(-85,100)(-85,100)(-90,90)

\qbezier(-85,100)(-85,100)(-70,70) \qbezier(-85,70)(-85,70)(-80,80)
\qbezier(-70,100)(-70,100)(-75,90)
\qbezier(-70,100)(-70,100)(-55,70) \qbezier(-70,70)(-70,70)(-65,80)
\qbezier(-55,100)(-55,100)(-60,90) \qbezier(-55,40)(-55,40)(-40,70)
\qbezier(-55,70)(-55,70)(-50,60) \qbezier(-40,40)(-40,40)(-45,50)
\qbezier(-40,40)(-40,40)(-25,70) \qbezier(-25,40)(-25,40)(-30,50)
\qbezier(-40,70)(-40,70)(-35,60)

\qbezier(-25,40)(-25,40)(-10,70) \qbezier(-10,40)(-10,40)(-15,50)
\qbezier(-25,70)(-25,70)(-20,60)

\put(-80,110){\makebox(0,0)[cc]{$\scriptstyle 2i+1$}}
\put(-30,30){\makebox(0,0)[cc]{$\scriptstyle 2j+1$}}
\end{picture}}
\put(120,0){\begin{picture}(0,130)(-40,0) \thicklines

\qbezier(-100,40)(-100,40)(-55,40) \qbezier(-25,40)(-25,40)(0,40)
\qbezier(-100,70)(-100,70)(-95,70) \qbezier(-65,70)(-65,70)(-55,70)
\qbezier(-25,70)(-25,70)(0,70) \qbezier(-100,100)(-100,100)(-95,100)
\qbezier(-65,100)(-65,100)(0,100)
\qbezier(-100,130)(-100,130)(0,130)

\put(0,85){\oval(40,30)[r]} \put(0,85){\oval(100,90)[r]}
\put(-100,55){\oval(40,30)[l]} \put(-100,115){\oval(40,30)[l]}

\put(-20,35){\oval(20,20)[b]} \put(-20,75){\oval(20,20)[t]}
\qbezier(-30,35)(-30,35)(-30,75) \qbezier(-10,45)(-10,45)(-10,65)

\put(-60,65){\oval(20,20)[b]} \put(-60,105){\oval(20,20)[t]}
\qbezier(-70,65)(-70,65)(-70,105) \qbezier(-50,75)(-50,75)(-50,95)

\qbezier(-95,100)(-95,100)(-80,70) \qbezier(-97,70)(-95,70)(-90,80)
\qbezier(-80,100)(-80,100)(-85,90)

\qbezier(-75,40)(-75,40)(-50,40) \qbezier(-25,40)(-25,40)(0,40)
\qbezier(-80,70)(-80,70)(-75,70) \qbezier(-65,70)(-65,70)(-50,70)
\qbezier(-80,100)(-80,100)(-75,100)

\qbezier(-50,40)(-50,40)(-35,70) \qbezier(-35,40)(-35,40)(-40,50)
\qbezier(-50,70)(-50,70)(-45,60)

\end{picture}}
\end{picture} \caption{$K(4ij+2i+2j+2,2j+1)$ and
$\mathcal B^{\prime \prime}$.} \label{Borromean:2prime:fig}
\end{center}
\end{figure}
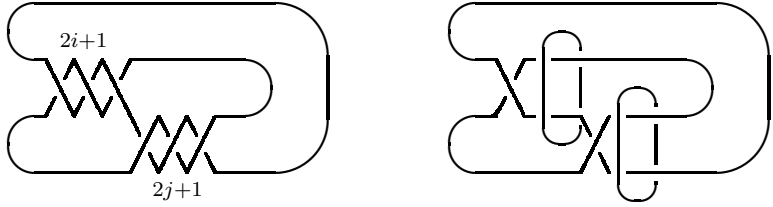

The proof is complete.\end{proof}

Propositions~\ref{obvious:prop},~\ref{limit:prop} and~\ref{l2:prop}
imply that
\begin{equation}
\lim_{k,m\to\infty} c(M_n(km+1,m)) \geqslant \left( 1 -
\frac{4\pi^2}{n^2} \right)^{3/2} \cdot
\frac{\vol(S^3\setminus\mathcal B)}{v_3}\cdot n
\end{equation}
with the factor multiplying $n$ in the right-hand side of the
inequality tending to $=7.21985\ldots$ as $n$ tends to $\infty$. But
fixing small $k$ and $m$, and using the computer program
SnapPea~\cite{snappea} to calculate the volume of $K(km+1,m)$, one
gets more specific lower bounds depending on $n$. For instance, with
notation again taken from~\cite{Rolfsen}, one can consider $K(5,2) =
4_1$, the figure-eight knot, $K(7,3)=5_2$, and $K(9,4)=6_1$, with
$\vol(S^3 \setminus 4_1) = 2 v_3$, $\vol(S^3 \setminus 5_2) =
2.81812\ldots$ , and $\vol(S^3 \setminus 6_1) = 3.16396 \ldots$ .
Propositions~\ref{obvious:prop} and \ref{limit:prop} imply for $n
\geqslant 7$ that
\begin{equation} \label{estimate:4_1}
c(M_n(5,2)) >
\left( 1 - \frac{4\pi^2}{n^2}\right)^{3/2}   \cdot 2n,
\end{equation}
(as a matter of fact, using an explicit formula for
$\vol(M_n(5,2))$, it was already shown in \cite{MaPV} that for
sufficiently large $n$ one has $c(M_n(5,2))> 2n$), and that
\begin{equation} \label{estimate:5_2}
c(M_n(7,3))>\left( 1 - \frac{4\pi^2}{n^2}\right)^{3/2} \cdot 2.77664 \ldots\ \cdot n,
\end{equation}
and for $n \geqslant 6$ that
\begin{equation} \label{estimate:6_1}
c(M_n(9,4))>\left( 1 - \frac{2\sqrt{2} \pi^2}{n^2} \right)^{3/2} \cdot
3.11739\ \ldots \cdot n.
\end{equation}

Turning to the upper estimate for the complexity of $M_n(p,q)$, we now note that in the special
case $\ell(p,q)=2$ the bound
$n(p-1)$ given by~(\ref{upper:estimate}) can also be significantly improved using a
more specific fundamental polyhedron instead of the Minkus polyhedron:

\begin{prop}\label{better:upper:prop}
Let $k,m\geqslant 2$ be integers. Suppose they are not both odd,
so $K(km+1,m)$ is a knot. Then, with the usual notation,
\begin{equation}
c(M_n(km+1,m))\leqslant n\cdot (\min\{k,m\}+k+m-3) \quad \forall
n.
\end{equation}
\end{prop}

\begin{proof}
It follows from~\cite{KK} that $M_n(mk+1,m)$ can be realized by
gluing together in pairs the faces of a polyhedron with $4n$ faces,
half being $(k+1)$-edged and half being $(m+1)$-edged polygons. More
precisely, this polyhedron is obtained by taking $n$ polygons with
$k+1$ (respectively, $m+1$) edges cyclically arranged around the
North (respectively, South) pole of the sphere, and $2n$ polygons
($n$ with $k+1$ and $n$ with $m+1$ edges) in the remaining
equatorial belt. In addition, each polygon incident to a pole is
glued to one in the equatorial belt.\footnote{As a minor fact we
note that there are misprints in Figg.~1 and~2 of~\cite{KK} for the
case where the integers involved have different parity, and in fact
the boundary patterns of $F_i$ and $\overline{F}_i$ do not match.
Using the notation of~\cite{KK}, so that the integers involved are
$m=2k+1$ and $s=2\ell$, one way of fixing these figures is as
follows. Keep calling $\ldots,F_i,F_{i+1},\ldots$ from left to right
the $m$-gons incident to the North pole $N$, so that $F_i$ has the
edges $x_i$ on its left and $x_{i+1}$ on its right, both emanating
from $N$. Similarly, call $\ldots,K_i,K_{i+1},\ldots$ from left to
right the $s$-gons incident to the South pole $S$, so that $K_i$ has
the edges $y_i$ on its left and $y_{i+1}$ on its right, both
emanating from $S$. Then the only $m$-gon adjacent to both $F_i$ and
$K_i$ should be $\overline{F}_{i+2}$, not $\overline{F}_i$, while
the only $s$-gon adjacent to both $K_i$ and $F_{i+1}$ should be
$\overline{K}_i$, as in~\cite{KK}. Now the boundary pattern of $F_i$
should be given, starting from $N$ and proceeding counterclockwise,
by the word
$$x_iy_{i-1}^{-1}x_{i+2k-1}^{-1}x_{i+2k-2}\cdots
x_{i+3}^{-1}x_{i+2}x_{i+1}^{-1},$$ while the boundary pattern for
$K_i$ should be given, starting from $S$ and proceeding clockwise,
by the word
$$y_ix_{i+2k-1}\Big(y_iy_{i+1}^{-1}\Big)^{\ell -1},$$
which allows to reconstruct the edge labelling completely.} Just
as in Lemma~3.1 of~\cite{MaPV}, we can now triangulate the
polygons incident to the poles by taking diagonals emanating from
the poles, and the polygons in the equatorial belt so that the
triangulations are matched under the gluing. If we now subdivide
the whole polyhedron by taking cones from the North pole, the
number of tetrahedra we obtain is given by the number of triangles
not incident to this pole, which is
$$n\cdot (k+1-2)+2n\cdot (m+1-2)=n\cdot (k+2m-3).$$
Similarly, if we take cones from the South pole we get $n\cdot (2k+m-3)$ tetrahedra,
and the conclusion readily follows.
\end{proof}

Note that~(\ref{upper:estimate}) gives $n\cdot km$ as an upper
estimate for $c(M_n(km+1,m))$, so the previous proposition always
gives a stronger bound. For instance for the knots $K(5,2) = 4_1$,
$K(7,3)=5_2$, and $K(9,4)=6_1$, namely those considered above
in~(\ref{estimate:4_1})--(\ref{estimate:6_1}),
formula~(\ref{upper:estimate}) gives as upper bounds $4n$, $6n$ and
$8n$ respectively, while Proposition~\ref{better:upper:prop} gives
$3n$, $4n$ and $5n$ respectively.

\end{document}